\documentclass[a4paper,11pt]{article}
\usepackage{graphicx}
\usepackage{fancybox}
\usepackage{amssymb}
\usepackage[all]{xy}
\usepackage{color}
\newtheorem{lemma}{Lemma}
\newtheorem{proposition}{Proposition}
\newtheorem{definition}{Definition}
\newtheorem{example}{Example}
\def\Rset{\mathrm{I\!R}}

\def\rien{\rule{0pt}{0pt}}
\def\DSmTchapter{paper}
\begin{document}
\begin{center}
\textbf{{\Huge Order in DSmT; definition of continuous DSm models}}\\
\vspace{1cm}
Fr\'ed\'eric Dambreville\\
D\'el\'egation G\'en\'erale pour l'Armement, DGA/CEP/GIP/SRO\\
\noindent
16 Bis, Avenue Prieur de la C\^ote d'Or\\
Arcueil, F 94114, France\\
Web: {\tt http://www.FredericDambreville.com}\\
Email: {\tt http://email.FredericDambreville.com}
\end{center}
\paragraph{Abstract.}
When implementing the DSmT, a difficulty may arise from the possible huge dimension of hyperpower sets, which are indeed free structures.
However, it is possible to reduce the dimension of these structures by involving logical constraints.
In this \DSmTchapter, the logical constraints will be related to a predefined order over the logical propositions.
The use of such orders and their resulting logical constraints will ensure a great reduction of the model complexity.
Such results will be applied to the definition of continuous DSm models.
In particular, a simplified description of the continuous impreciseness is considered, based on  impreciseness intervals of the sensors.
From this viewpoint, it is possible to manage the contradictions between continuous sensors in a DSmT manner, while the complexity of the model stays handlable.
\section{Introduction}
Recent advances\cite{DSmTBook1} in the Dezert Smarandache Theory have shown that this theory was able to handle the contradiction between propositions in a quite flexible way.
This new theory has been already applied in different domains; \emph{e.g.}:
\begin{itemize}
\item Data asociation in target tracking \cite{tchamova}\,,
\item Environmental prediction \cite{corgne}\,.
\end{itemize}
Although free DSm models are defined over hyperpower sets, which sizes evolve exponentially with the number of \emph{atomic} propositions, it appears that the manipulation of the fusion rule is still manageable for practical problems reasonnably well shaped.
Moreover, the hybrid DSm models are of lesser complexity.
\\[5pt]
If DSmT works well for discret spaces, the manipulation of continous DSm models is still an unknown.
Nevertheless, the management of continuous data is an issue of main interest.
It is necessary for implementing a true fusion engine for localization informations;
and associated with a principle of conditionning, it will be a main ingredient for implementing filters for the localization.
But a question first arises: \emph{what could be an hyperpower set for a continuous DSm model?}
Such first issue does not arises so dramatically in Dempster Shafer Theory or for Transfer Belief Models\cite{smetskennes}.
In DST, a continuous proposition could just be a measurable subset.
On the other hand, a free DSm model, defined over an hyperpower set, will imply that any pair of propositions will have a non empty intersection.
This is desappointing, since the notion of \emph{point} (a minimal non empty proposition) does not exist anymore in an hyperpower set.
\\[5pt]
But even if it is possible to define a continuous propositional model in DST/TBM, the manipulation of continuous basic belief assignment is still an issue\cite{risticsmets}\cite{strat}.
In \cite{risticsmets}\,, Ristic and Smets proposed a restriction of the bba to intervals of $\Rset$\,.
It was then possible to derive a mathematical relation between a continuous bba density and its $\mathrm{Bel}$ function.
\\[5pt]
In this \DSmTchapter, the construction of continuous DSm models is proposed.
This construction is based on a constrained model, where the logical contraints are implied by the definition of an order relation over the propositions.
\\[3pt]
A one-dimension DSm model will be implemented, where the definition of the basic belief assignment relies on a \emph{generalized notion of intervals}.
Although this construction has been fulfilled on a different ground, it shares some surprizing similarities with Ristic and Smets viewpoint.
As in \cite{risticsmets}, the bba will be seen as density defined over a 2-dimension measurable space.
We will be able to derive the Belief function from the basic belief assignment, by applying an integral computation.
At last, the conjunctive fusion operator, $\oplus$, is derived by a rather simple integral computation.
\\[9pt]
Section~\ref{f2k5:sectionIntroDSmT} makes a quick introduction of the Dezert Smarandache Theory.
Section~\ref{f2k5:sectionTowardOrderedModel} is about ordered DSm models.
In section~\ref{f2k5:sectionTowardContinuousModel}, a continuous DSm model is defined.
This method is restricted to only one dimension.
The related computation methods are detailed.
In section~\ref{f2k5:sectionImplementation}, our algorithmic implementation is described and an example of computation is given.
%
%
The paper is then concluded.
\section{A short introduction to the DSmT}
\label{f2k5:sectionIntroDSmT}
The theory and its meaning are widely explained in \cite{DSmTBook1}.
However, we will particularly focus on the notion of hyperpower sets, since this notion is fundamental subsequently.
\\[5pt]
The \emph{Dezert Smarandache Theory} belongs to the family of \emph{Evidence Theories}.
As the \emph{Dempster Shafer Theory}\cite{dempster}\cite{shafer} or the \emph{Transferable Belief Models}\cite{smetskennes}, the DSmT is a framework for fusing belief informations, originating from independent sensors.
However, free DSm models are defined over Hyperpower sets, which are \emph{fully open-world extensions} of sets.
It is possible to restrict this full open-world hypothesis by adding propositional constraints, resulting in the definition of an \emph{hybrid Dezert Smarandache model}.
\\[5pt]
The notion of hyperpower set is thus a fundamental ingredient of the DSmT.
Hyperpower sets could be considered as a free pre-Boolean algebra.
As these structures will be of main importance subsequently, the next sections are devoted to introduce them in details.
As a prerequisite, the notion of Boolean algebra is quickly introduced now.
\subsection{Boolean algebra}
\paragraph{Definition.}
A Boolean algebra is a sextuple $(\Phi,\wedge,\vee,\neg,\bot,\top)$ such that:
\begin{itemize}
\item $\Phi$ is a set, called set of propositions,
\item $\bot,\top$ are specific propositions of $\Phi$, respectively called \emph{false} and \emph{true},
\item $\neg:\Phi\rightarrow \Phi$ is a unary operator,
\item $\wedge:\Phi\times \Phi\rightarrow \Phi$ and $\vee:\Phi\times \Phi\rightarrow \Phi$ are binary operators,
\end{itemize}
and verifying the following properties:
\begin{description}
\item[A1.] $\wedge$ and $\vee$ are commutative:
\vspace{-5pt}
$$
\rien\hspace{-25pt}
\forall \phi,\psi\in\Phi,\ 
\phi\wedge\psi=\psi\wedge\phi
\mbox{ and }
\phi\vee\psi=\psi\vee\phi\;,
$$
\item[A2.] $\wedge$ and $\vee$ are associative:
\vspace{-5pt}
$$
\rien\hspace{-25pt}
\forall \phi,\psi,\eta\in\Phi,\ 
(\phi\wedge\psi)\wedge\eta=\phi\wedge(\psi\wedge\eta)
\mbox{ and }
(\phi\vee\psi)\vee\eta=\phi\vee(\psi\vee\eta)\;,
$$
\item[A3.] $\top$ is neutral for $\wedge$ and $\bot$ is neutral for $\vee$:
\vspace{-5pt}
$$
\rien\hspace{-25pt}
\forall \phi\in\Phi,\ 
\phi\wedge\top=\phi
\mbox{ and }
\phi\vee\bot=\phi\;,
$$
\item[A4.] $\wedge$ and $\vee$ are distributive for each other:
\vspace{-5pt}
$$
\rien\hspace{-25pt}
\forall \phi,\psi,\eta\in\Phi,\ 
\phi\wedge(\psi\vee\eta)=(\phi\wedge\psi)\vee(\phi\wedge\eta)
\mbox{ and }
\phi\vee(\psi\wedge\eta)=(\phi\vee\psi)\wedge(\phi\vee\eta)
\;,
$$
\item[A5.] $\neg$ defines the complement of any proposition:
\vspace{-5pt}
$$
\rien\hspace{-25pt}
\forall \phi\in\Phi,\ 
\phi\wedge\neg\phi=\bot
\mbox{ and }
\phi\vee\neg\phi=\top\;.
$$
\end{description}
The Boolean algebra $(\Phi,\wedge,\vee,\neg,\bot,\top)$ will be also refered to as the \emph{Boolean algebra $\Phi$}, the strucure being thus implied.
An order relation $\subset$ is defined over $\Phi$ by:
$$
\forall\phi,\psi\in\Phi\,,\quad \phi\subset\psi \stackrel{\Delta}{\iff}\phi\wedge\psi=\phi\;.
$$
\paragraph{Fundamental examples.}
The following examples are two main conceptions of Boolean algebra.
\begin{example}
Let $\Omega$ be a set and $\mathcal{P}(\Omega)$ be the set of its subsets.
For any $A\subset\Omega$, denote $\sim A=\Omega\setminus A$ its complement.
Then $\bigl(\mathcal{P}(\Omega),\cap,\cup,\sim,\emptyset,\Omega\bigr)$ is a Boolean algebra.
\end{example}
The proof is immediate by verifying the properties A1 to A5.
\begin{example}
For any $i\in\{1,\dots,n\}$, let $\theta_i=\{0,1\}^{i-1}\times\{0\}\times\{0,1\}^{n-i}$\,.
Let $\Theta=\{\theta_1,\dots,\theta_n\}$ and denote $\bot=\emptyset$, $\top=\{0,1\}^n$ and $\mathcal{B}(\Theta)=\mathcal{P}\bigl(\{0,1\}^n\bigr)$\,.
Define the operators $\wedge$, $\vee$ and $\neg$ by $\phi\wedge\psi=\phi\cap\psi$, $\phi\vee\psi=\phi\cup\psi$ and $\neg\phi=\top\setminus\phi$ for any $\phi,\psi\in\mathcal{B}(\Theta)$\,.
Then $\bigl(\mathcal{B}(\Theta),\wedge,\vee,\neg,\bot,\top\bigr)$ is a Boolean algebra.
\end{example}
\begin{figure}[h!]
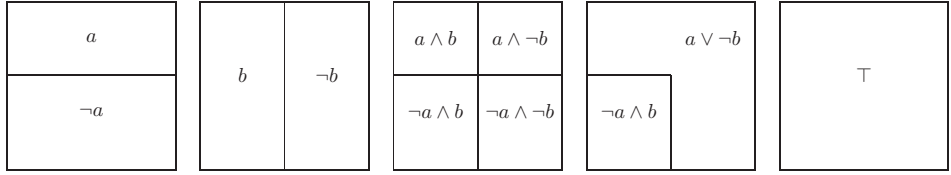

\centering
\scalebox{.70}{\begin{tabular}{@{}c@{ \ }c@{ \ }c@{ \ }c@{ \ }c@{}}
\xy
(0,0)*+{\rien},
(0,0)="box00",
(16,0)="box10",
(32,0)="box20",
(0,-13.9)="box01",
(16,-13.9)="box11",
(32,-13.9)="box21",
(0,-32)="box02",
(16,-32)="box12",
(32,-32)="box22",
\ar @{-}"box00";"box20",
\ar @{-}"box01";"box21",
\ar @{-}"box02";"box22",
\ar @{-}"box00";"box02",
\ar @{-}"box20";"box22",
\POS "box22",
(16,-6.9)*+{a},
(16,-20.8)*+{\neg a},
(32,-32)*+{\rien}
\endxy
&
\xy
(0,0)*+{\rien},
(0,0)="box00",
(16,0)="box10",
(32,0)="box20",
(0,-13.9)="box01",
(16,-13.9)="box11",
(32,-13.9)="box21",
(0,-32)="box02",
(16,-32)="box12",
(32,-32)="box22",
\ar @{-}"box00";"box20",
\ar @{-}"box02";"box22",
\ar @{-}"box00";"box02",
\ar @{-}"box10";"box12",
\ar @{-}"box20";"box22",
\POS "box22",
(8,-13.9)*+{b},
(24,-13.9)*+{\neg b},
(32,-32)*+{\rien}
\endxy
&
\xy
(0,0)*+{\rien},
(0,0)="box00",
(16,0)="box10",
(32,0)="box20",
(0,-13.9)="box01",
(16,-13.9)="box11",
(32,-13.9)="box21",
(0,-32)="box02",
(16,-32)="box12",
(32,-32)="box22",
\ar @{-}"box00";"box20",
\ar @{-}"box01";"box21",
\ar @{-}"box02";"box22",
\ar @{-}"box00";"box02",
\ar @{-}"box10";"box12",
\ar @{-}"box20";"box22",
\POS "box22",
(8,-6.9)*+{a\wedge b},
(24,-6.9)*+{a\wedge\neg b},
(8,-20.8)*+{\neg a\wedge b},
(24,-20.8)*+{\neg a\wedge\neg b},
(32,-32)*+{\rien}
\endxy
&
\xy
(0,0)*+{\rien},
(0,0)="box00",
(16,0)="box10",
(32,0)="box20",
(0,-13.9)="box01",
(16,-13.9)="box11",
(32,-13.9)="box21",
(0,-32)="box02",
(16,-32)="box12",
(32,-32)="box22",
\ar @{-}"box00";"box20",
\ar @{-}"box01";"box11",
\ar @{-}"box02";"box22",
\ar @{-}"box00";"box02",
\ar @{-}"box11";"box12",
\ar @{-}"box20";"box22",
\POS "box22",
(24,-6.9)*+{a\vee\neg b},
(8,-20.8)*+{\neg a\wedge b},
(32,-32)*+{\rien}
\endxy
&
\xy
(0,0)*+{\rien},
(0,0)="box00",
(16,0)="box10",
(32,0)="box20",
(0,-13.9)="box01",
(16,-13.9)="box11",
(32,-13.9)="box21",
(0,-32)="box02",
(16,-32)="box12",
(32,-32)="box22",
\ar @{-}"box00";"box20",
\ar @{-}"box02";"box22",
\ar @{-}"box00";"box02",
\ar @{-}"box20";"box22",
\POS "box22",
(16,-13.9)*+{\top},
(32,-32)*+{\rien}
\endxy
\end{tabular}}
\caption{Boolean algebra $\mathcal{B}\bigl(\{a,b\}\bigr)$; (partial)}
\label{DSmTBook2:DMB:AppendBalg:1}
\end{figure}
The second example seems just like a rewriting of the first one, but it is of the most importance.
It is called the \emph{free Boolean algebra generated by the set of atomic propositions $\Theta$}.
Figure~\ref{DSmTBook2:DMB:AppendBalg:1} shows the structure of such algebra, when $n=2$.
The free Boolean algebra $\mathcal{B}(\Theta)$ is deeply related to the classical propositional logic: it gives the (logical) equivalence classes of the propositions generated from the atomic propositions of $\Theta$. 
Although we give here an explicit definition of $\mathcal{B}(\Theta)$ by means of its binary coding $\mathcal{P}\bigl(\{0,1\}^n\bigr)$, the truely rigorous definition of $\mathcal{B}(\Theta)$ is made by means of the logical equivalence (which is out of the scope of this presentation).
Thus, the binary coding of the atomic propositions $\theta_i\in\Theta$ is only implied.
\paragraph{Fondamental proposition.}
\begin{proposition}
Any Boolean algebra is isomorph to a Boolean algebra derived from a set, \emph{i.e.} $\bigl(\mathcal{P}(\Omega),\cap,\cup,\sim,\emptyset,\Omega\bigr)$\,.
\end{proposition}
Proofs should be found in any good reference; see also \cite{DSmTbook2:DMB:boolAlgebra}.
\subsection{Hyperpower sets}
\paragraph{Definition of hyperpower set.}
Let $\Theta$ be a finite set of atomic propositions, and $\bigl(\mathcal{B}(\Theta),\wedge,\vee,\neg,\bot,\top\bigr)$ be the free Boolean algebra generated by $\Theta$.
For any $\Sigma\subset\mathcal{P}(\Theta)$, define $\varphi(\Sigma)$, element of $\mathcal{B}(\Theta)$, by \mbox{$\varphi(\Sigma)=\bigvee_{\sigma\in \Sigma}\bigwedge_{\theta\in \sigma}\theta\;.$}\footnote{It is assumed $\bigvee_{\phi\in\emptyset}=\bot$ and $\bigwedge_{\phi\in\emptyset}=\top$\,.}
The set $<\Theta>=\bigl\{\varphi(\Sigma)\;\big/\;\Sigma\subset\mathcal{P}(\Theta)\bigr\}$ is called hyperpower set generated by~$\Theta$.
\\[5pt]
It is noticed that both $\bot=\varphi(\emptyset)$ and $\top=\varphi\bigl(\mathcal{P}(\Theta)\bigr)$ are elements of $<\Theta>$.
Figure~\ref{DSmTBook2:DMB:AppendBalg:2} shows the structure of the hyperpower set, when $n=2$.
Typically, it appears that the elements of the hyperpower set are built only from $\neg$-free components.
\begin{figure}[h!]
\centering
\scalebox{1}{\begin{tabular}{@{}c@{ \ }c@{ \ }c@{ \ }c@{ \ }c@{ \ }c@{}}
\xy
(0,0)*+{\rien},
(0,0)="box00",
(8,0)="box10",
(16,0)="box20",
(0,-8)="box01",
(8,-8)="box11",
(16,-8)="box21",
(0,-16)="box02",
(8,-16)="box12",
(16,-16)="box22",
\ar @{--}"box00";"box20",
\ar @{--}"box02";"box22",
\ar @{--}"box00";"box02",
\ar @{--}"box20";"box22",
\POS "box22",
%
(16,-16)*+{\rien}
\endxy
&
\xy
(0,0)*+{\rien},
(0,0)="box00",
(8,0)="box10",
(16,0)="box20",
(0,-8)="box01",
(8,-8)="box11",
(16,-8)="box21",
(0,-16)="box02",
(8,-16)="box12",
(16,-16)="box22",
\ar @{-}"box00";"box10",
\ar @{--}"box10";"box20",
\ar @{-}"box01";"box11",
\ar @{--}"box02";"box22",
\ar @{-}"box00";"box01",
\ar @{--}"box01";"box02",
\ar @{-}"box10";"box11",
\ar @{--}"box20";"box22",
\POS "box22",
%
(16,-16)*+{\rien}
\endxy
&
\xy
(0,0)*+{\rien},
(0,0)="box00",
(8,0)="box10",
(16,0)="box20",
(0,-8)="box01",
(8,-8)="box11",
(16,-8)="box21",
(0,-16)="box02",
(8,-16)="box12",
(16,-16)="box22",
\ar @{-}"box00";"box20",
\ar @{-}"box01";"box21",
\ar @{--}"box02";"box22",
\ar @{-}"box00";"box01",
\ar @{--}"box01";"box02",
\ar @{-}"box20";"box21",
\ar @{--}"box21";"box22",
\POS "box22",
%
(16,-16)*+{\rien}
\endxy
&
\xy
(0,0)*+{\rien},
(0,0)="box00",
(8,0)="box10",
(16,0)="box20",
(0,-8)="box01",
(8,-8)="box11",
(16,-8)="box21",
(0,-16)="box02",
(8,-16)="box12",
(16,-16)="box22",
\ar @{-}"box00";"box10",
\ar @{--}"box10";"box20",
\ar @{-}"box02";"box12",
\ar @{--}"box12";"box22",
\ar @{-}"box00";"box02",
\ar @{-}"box10";"box12",
\ar @{--}"box20";"box22",
\POS "box22",
%
(16,-16)*+{\rien}
\endxy
&
\xy
(0,0)*+{\rien},
(0,0)="box00",
(8,0)="box10",
(16,0)="box20",
(0,-8)="box01",
(8,-8)="box11",
(16,-8)="box21",
(0,-16)="box02",
(8,-16)="box12",
(16,-16)="box22",
\ar @{-}"box00";"box20",
\ar @{-}"box20";"box21",
\ar @{--}"box21";"box22",
\ar @{-}"box21";"box11",
\ar @{-}"box00";"box02",
\ar @{-}"box02";"box12",
\ar @{--}"box12";"box22",
\ar @{-}"box12";"box11",
\POS "box22",
%
(16,-16)*+{\rien}
\endxy
&
\xy
(0,0)*+{\rien},
(0,0)="box00",
(8,0)="box10",
(16,0)="box20",
(0,-8)="box01",
(8,-8)="box11",
(16,-8)="box21",
(0,-16)="box02",
(8,-16)="box12",
(16,-16)="box22",
\ar @{-}"box00";"box20",
\ar @{-}"box02";"box22",
\ar @{-}"box00";"box02",
\ar @{-}"box20";"box22",
\POS "box22",
%
(16,-16)*+{\rien}
\endxy
\\
$\bot$&$a\wedge b$&$a$&$b$&$a\vee b$&$\top$
\end{tabular}}
\caption{Hyperpower set $<a,b>=\{\bot, a\wedge b, a, b,a\vee b, \top\}$}
\label{DSmTBook2:DMB:AppendBalg:2}
\end{figure}
\begin{example}\label{DSMTb2:dmb:ex2}
Hyperpower set generated by $\Theta=\{a,b,c\}$.
$$\begin{array}{@{}l@{}}\displaystyle
<a,b,c>=\bigl\{\bot,a,b,c,
a\wedge b\wedge c,
a\wedge b, b\wedge c, c\wedge a,a\vee b\vee c,\top
\\\qquad\qquad\qquad
a\vee b, b\vee c, c\vee a,
(a\wedge b)\vee c, (b\wedge c)\vee a, (c\wedge a)\vee b,
\\\qquad\qquad\qquad
(a\vee b)\wedge c, (b\vee c)\wedge a, (c\vee a)\wedge b,
(a\wedge b)\vee(b\wedge c)\vee(c\wedge a)
\bigr\}
\end{array}$$
The following table associates some $\Sigma\subset\mathcal{P}(\Theta)$ to their related hyperpower element $\varphi(\Sigma)$.
This table is partial; there is indeed $256$ possible choices for $\Sigma$.
It appears that $\varphi$ is not one-to-one:
{\small$$
\begin{array}{c|c|c}
\Sigma &\varphi(\Sigma)&\mbox{reduced form in }<\Theta>
\\\hline\hline
\emptyset
&
\bot
&
\bot
\\\hline
\{\emptyset\}
&
\top
&
\top
\\\hline
\bigl\{
\{a\}
;
\{b\}
;
\{c\}
\bigr\}
&
a\vee b\vee c
&
a\vee b\vee c
\\\hline
\bigl\{
\{a,b\}
;
\{b,c\}
;
\{c,a\}
\bigr\}
&
(a\wedge b)\vee(b\wedge c)\vee (c\wedge a)
&
(a\wedge b)\vee(b\wedge c)\vee (c\wedge a)
\\\hline
\bigl\{
\{a,c\}
;
\{b,c\}
;
\{a,b,c\}
\bigr\}
&
(a\wedge c)\vee(b\wedge c)\vee (a\wedge b\wedge c)
&
(a\vee b)\wedge c
\\\hline
\bigl\{
\{a,c\}
;
\{b,c\}
\bigr\}
&
(a\wedge c)\vee(b\wedge c)
&
(a\vee b)\wedge c
\end{array}
$$}
\end{example}
\emph{Remark.} In the DSmT book~1 \cite{DSmTBook1}, the hyperpower sets have been defined by means of the Smarandache encoding.
Our definition is quite related to this encoding.
In fact this encoding is just implied in the definition of $\varphi$.
\paragraph{Hyperpower set as a free pre-Boolean algebra.}
It is easy to verify on example~\ref{DSMTb2:dmb:ex2} that $<\Theta>$ is left unchanged by any application of the operators $\wedge$ and $\vee$.
For example:
$$
(a\wedge b)\wedge\bigl((b\wedge c)\vee a\bigr)=(a\wedge b\wedge b\wedge c)\vee(a\wedge b\wedge a)=a\wedge b\;.
$$
This result is formalized by the following proposition.
\begin{proposition}
\label{DSmTb2:dmb:Prop1}
Let $\phi,\psi\in<\Theta>$.
Then $\phi\wedge\psi\in<\Theta>$ and $\phi\vee\psi\in<\Theta>$.
\end{proposition}
\begin{description}
\item[Proof.]
Let $\phi,\psi\in<\Theta>$.\\
There are $\Sigma\subset\mathcal{P}(\Theta)$ and $\Gamma\subset\mathcal{P}(\Theta)$ such that $\phi=\varphi(\Sigma)$ and $\psi=\varphi(\Gamma)$\,.\\
By applying the definition of $\varphi$, it comes immediately:
$$
\varphi(\Sigma)\vee\varphi(\Gamma)=\bigvee_{\sigma\in \Sigma\cup\Gamma}\bigwedge_{\theta\in \sigma}\theta\;.
$$
It is also deduced:
$$
\varphi(\Sigma)\wedge\varphi(\Gamma)=\Biggl(\bigvee_{\sigma\in \Sigma}\bigwedge_{\theta\in \sigma}\theta\Biggr)
\wedge
\Biggl(\bigvee_{\gamma\in\Gamma}\bigwedge_{\theta\in \gamma}\theta\Biggr)
\;.
$$
By applying the distributivity, it comes:
$$
\varphi(\Sigma)\wedge\varphi(\Gamma)=\bigvee_{\sigma\in \Sigma}\bigvee_{\gamma\in\Gamma}\Biggl(
\biggl(\bigwedge_{\theta\in \sigma}\theta\biggr)
\wedge
\biggl(\bigwedge_{\theta\in \gamma}\theta\biggr)
\Biggr)
=
\bigvee_{(\sigma,\gamma)\in \Sigma\times\Gamma}\bigwedge_{\theta\in \sigma\cup\gamma}\theta
\;.
$$
Then $\varphi(\Sigma)\wedge\varphi(\Gamma)=\varphi(\Lambda)$\,, with $\Lambda=\bigl\{\sigma\cup\gamma\;\big/\;(\sigma,\gamma)\in \Sigma\times\Gamma\bigr\}$\,. 
\item[$\Box\Box\Box$]\rien
\end{description}
\emph{Corollary and definition.}
Proposition~\ref{DSmTb2:dmb:Prop1} implies that $\wedge$ and $\vee$ infer inner operations within $<\Theta>$\,.
As a consequence, $\bigl(<\Theta>,\wedge,\vee,\bot,\top\bigr)$ is an algebraic structure by itself.
Since it does not contains the negation~$\neg$, this structure is called the \emph{free pre-Boolean algebra generated by $\Theta$.}
\subsection{Pre-Boolean algebra}
\paragraph{Generality.}
Typically, a free algebra is an algebra where the only constraints are the intrinsic contraints which characterize its fundamental structures.
For example in a free Boolean algebra, the only constraints are A1 to A5, and there are no other constraints put on the propositions.
But conversely, it is indeed possible to derive any algebra by constraining its free counterpart.
This will be our approach for defining pre-Boolean algebra in general:
a pre-Boolean algebra will be a \emph{constrained} free pre-Boolean algebra.
Constraining a structure is a quite intuitive notion.
However, a precise mathematical definition needs the abstract notion of equivalence relations and classes.
Let us start with the intuition by introducing an example.
\begin{example}\label{DSMTb2:dmb:ex3}
Pre-Boolean algebra generated by $\Theta=\{a,b,c\}$ and constrained by $a\wedge b=a\wedge c$ and $a\wedge c=b\wedge c$.
\\
\emph{For coherence with forthcoming notations, these constraints will be designated by using the set of propositional pairs \mbox{$\Gamma=\bigl\{(a\wedge b,a\wedge c),(a\wedge c,b\wedge c)\bigr\}\,.$}}
\\
The idea is to start from the free pre-Boolean algebra $<a,b,c>$, propagate the constraints, and then reduce the propositions identified by the constraints.
\\[4pt]
It is first deduced $a\wedge b=a\wedge c=b\wedge c=a\wedge b\wedge c$.\\
It follows $(a\wedge b)\vee c=c$, $(b\wedge c)\vee a=a$ and $(c\wedge a)\vee b=b$.\\
Also holds $(a\vee b)\wedge c= (b\vee c)\wedge a= (c\vee a)\wedge b=
(a\wedge b)\vee(b\wedge c)\vee(c\wedge a)=a\wedge b\wedge c$\,.\\
By discarding these cases from the free structure $<a,b,c>$, it comes the following constrained pre-Boolean algebra:
$$
<a,b,c>_\Gamma=\bigl\{\bot,a\wedge b\wedge c,a,b,c,a\vee b, b\vee c, c\vee a,a\vee b\vee c,\top\bigr\}
$$
Of course, it is necessary to show that there is actually no further reduction in $<a,b,c>_\Gamma$.
This is done by expliciting a model; for example the structure of figure~\ref{dsmtb2:dmb:fig3}.
\begin{figure}[ht!]
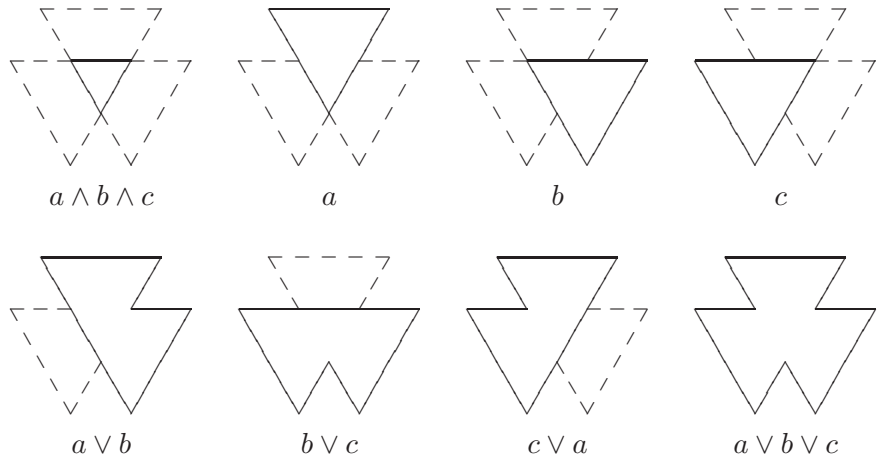

\centering
\scalebox{1}{\begin{tabular}{@{}cccc@{}}
\xy
(0,0)*+{\rien},
(4,0)="a1",
(20,0)="b1",
(12,-13.9)="c1",
(24,-6.9)="a2",
(16,-20.8)="b2",
(8,-6.9)="c2",
(8,-20.8)="a3",
(0,-6.9)="b3",
(16,-6.9)="c3",
\ar @{--}"a1";"b1",
\ar @{--}"b1";"c3",
\ar @{-}"c3";"c2",
\ar @{--}"c2";"a1",
\ar @{--}"a2";"b2",
\ar @{--}"b2";"c1",
\ar @{-}"c1";"c3",
\ar @{--}"c3";"a2",
\ar @{--}"a3";"b3",
\ar @{--}"b3";"c2",
\ar @{-}"c2";"c1",
\ar @{--}"c1";"a3",
\POS "c3",
(24,-20.8)*+{\rien}
\endxy
&
\xy
(0,0)*+{\rien},
(4,0)="a1",
(20,0)="b1",
(12,-13.9)="c1",
(24,-6.9)="a2",
(16,-20.8)="b2",
(8,-6.9)="c2",
(8,-20.8)="a3",
(0,-6.9)="b3",
(16,-6.9)="c3",
\ar @{-}"a1";"b1",
\ar @{-}"b1";"c3",
\ar @{-}"c2";"a1",
\ar @{--}"a2";"b2",
\ar @{--}"b2";"c1",
\ar @{-}"c1";"c3",
\ar @{--}"c3";"a2",
\ar @{--}"a3";"b3",
\ar @{--}"b3";"c2",
\ar @{-}"c2";"c1",
\ar @{--}"c1";"a3",
\POS "c3",
(24,-20.8)*+{\rien}
\endxy
&
\xy
(0,0)*+{\rien},
(4,0)="a1",
(20,0)="b1",
(12,-13.9)="c1",
(24,-6.9)="a2",
(16,-20.8)="b2",
(8,-6.9)="c2",
(8,-20.8)="a3",
(0,-6.9)="b3",
(16,-6.9)="c3",
\ar @{--}"a1";"b1",
\ar @{--}"b1";"c3",
\ar @{-}"c3";"c2",
\ar @{--}"c2";"a1",
\ar @{-}"a2";"b2",
\ar @{-}"b2";"c1",
\ar @{-}"c3";"a2",
\ar @{--}"a3";"b3",
\ar @{--}"b3";"c2",
\ar @{-}"c2";"c1",
\ar @{--}"c1";"a3",
\POS "c3",
(24,-20.8)*+{\rien}
\endxy
&
\xy
(0,0)*+{\rien},
(4,0)="a1",
(20,0)="b1",
(12,-13.9)="c1",
(24,-6.9)="a2",
(16,-20.8)="b2",
(8,-6.9)="c2",
(8,-20.8)="a3",
(0,-6.9)="b3",
(16,-6.9)="c3",
\ar @{--}"a1";"b1",
\ar @{--}"b1";"c3",
\ar @{-}"c3";"c2",
\ar @{--}"c2";"a1",
\ar @{--}"a2";"b2",
\ar @{--}"b2";"c1",
\ar @{-}"c1";"c3",
\ar @{--}"c3";"a2",
\ar @{-}"a3";"b3",
\ar @{-}"b3";"c2",
\ar @{-}"c1";"a3",
\POS "c3",
(24,-20.8)*+{\rien}
\endxy
\\
$a\wedge b\wedge c$&$a$&$b$&$c$
\vspace{3pt}
\\
\xy
(0,0)*+{\rien},
(4,0)="a1",
(20,0)="b1",
(12,-13.9)="c1",
(24,-6.9)="a2",
(16,-20.8)="b2",
(8,-6.9)="c2",
(8,-20.8)="a3",
(0,-6.9)="b3",
(16,-6.9)="c3",
\ar @{-}"a1";"b1",
\ar @{-}"b1";"c3",
\ar @{-}"c2";"a1",
\ar @{-}"a2";"b2",
\ar @{-}"b2";"c1",
\ar @{-}"c3";"a2",
\ar @{--}"a3";"b3",
\ar @{--}"b3";"c2",
\ar @{-}"c2";"c1",
\ar @{--}"c1";"a3",
\POS "c3",
(24,-20.8)*+{\rien}
\endxy
&
\xy
(0,0)*+{\rien},
(4,0)="a1",
(20,0)="b1",
(12,-13.9)="c1",
(24,-6.9)="a2",
(16,-20.8)="b2",
(8,-6.9)="c2",
(8,-20.8)="a3",
(0,-6.9)="b3",
(16,-6.9)="c3",
\ar @{--}"a1";"b1",
\ar @{--}"b1";"c3",
\ar @{-}"c3";"c2",
\ar @{--}"c2";"a1",
\ar @{-}"a2";"b2",
\ar @{-}"b2";"c1",
\ar @{-}"c3";"a2",
\ar @{-}"a3";"b3",
\ar @{-}"b3";"c2",
\ar @{-}"c1";"a3",
\POS "c3",
(24,-20.8)*+{\rien}
\endxy
&
\xy
(0,0)*+{\rien},
(4,0)="a1",
(20,0)="b1",
(12,-13.9)="c1",
(24,-6.9)="a2",
(16,-20.8)="b2",
(8,-6.9)="c2",
(8,-20.8)="a3",
(0,-6.9)="b3",
(16,-6.9)="c3",
\ar @{-}"a1";"b1",
\ar @{-}"b1";"c3",
\ar @{-}"c2";"a1",
\ar @{--}"a2";"b2",
\ar @{--}"b2";"c1",
\ar @{-}"c1";"c3",
\ar @{--}"c3";"a2",
\ar @{-}"a3";"b3",
\ar @{-}"b3";"c2",
\ar @{-}"c1";"a3",
\POS "c3",
(24,-20.8)*+{\rien}
\endxy
&
\xy
(0,0)*+{\rien},
(4,0)="a1",
(20,0)="b1",
(12,-13.9)="c1",
(24,-6.9)="a2",
(16,-20.8)="b2",
(8,-6.9)="c2",
(8,-20.8)="a3",
(0,-6.9)="b3",
(16,-6.9)="c3",
\ar @{-}"a1";"b1",
\ar @{-}"b1";"c3",
\ar @{-}"c2";"a1",
\ar @{-}"a2";"b2",
\ar @{-}"b2";"c1",
\ar @{-}"c3";"a2",
\ar @{-}"a3";"b3",
\ar @{-}"b3";"c2",
\ar @{-}"c1";"a3",
\POS "c3",
(24,-20.8)*+{\rien}
\endxy
\\
$a\vee b$&$b\vee c$&$c\vee a$&$a\vee b\vee c$
\end{tabular}}
\caption{Pre-Boolean algebra $<a,b,c>_\Gamma$; ($\bot$ and $\top$ are omitted)}
\label{dsmtb2:dmb:fig3}
\end{figure}
\end{example}
For the reader not familiar with the notion of equivalence classes, the following construction is just a mathematical formalization of the contraint propagation which has been described in example~\ref{DSMTb2:dmb:ex3}.
Now, it is first introduced the notion of morphism between structures.
\paragraph{Magma.}
A \emph{$(\wedge,\vee,\bot,\top)$-magma}, also called \emph{magma} for short, is a quintuple $(\Phi,\wedge,\vee,\bot,\top)$ where $\Phi$ is a set of propositions, $\wedge$ and $\vee$ are binary operators on $\Phi$, and $\bot$ and $\top$ are two elements of $\Phi$.
\\[5pt]
The magma $(\Phi,\wedge,\vee,\bot,\top)$ may also be refered to as the magma $\Phi$, the structure being thus implied.
Notice that an hyperpower set is a magma.
\paragraph{Morphism.}
Let $(\Phi,\wedge,\vee,\bot,\top)$ and  $(\Psi,\wedge,\vee,\bot,\top)$ be two magma.
A morphism $\mu$ from the magma $\Phi$ to the magma $\Psi$ is a mapping from $\Phi$ to $\Psi$ such that:
\begin{itemize}
\item $\mu(\phi\wedge\psi)=\mu(\phi)\wedge\mu(\psi)$ and $\mu(\phi\vee\psi)=\mu(\phi)\vee\mu(\psi)$\,,
\item $\mu(\bot)=\bot$ and $\mu(\top)=\top$\,.
\end{itemize}
A morphism is an isomorphism if it is a bijective mapping.
In such case, the magma $\Phi$ and the magma $\Psi$ are said to be isomorph, which means that they share the same structure.
\\[5pt]
\emph{The notions of $(\wedge,\vee)$-magma and of $(\wedge,\vee)$-morphism are defined similarly by discarding $\bot$ and $\top$.}
\paragraph{Propagation relation.}
Let $<\Theta>$ be a free pre-Boolean algebra.
Let $\Gamma\subset<\Theta>\times<\Theta>$ be a set of propositional pairs;
for any pair $(\phi,\psi)\in\Gamma$ is defined the constraint $\phi=\psi$\,.
The propagation relation associated to the constraints, and also denoted $\Gamma$, is defined recursively by:
\begin{itemize}
\item $\phi\Gamma\phi$, for any $\phi\in<\Theta>$, 
\item If $(\phi,\psi)\in\Gamma$, then $\phi\Gamma\psi$ and $\psi\Gamma\phi$,
\item If $\phi\Gamma\psi$ and $\psi\Gamma\eta$, then $\phi\Gamma\eta$,
\item If $\phi\Gamma\eta$ and $\psi\Gamma\zeta$, then $(\phi\wedge\psi)\Gamma(\eta\wedge\zeta)$ and $(\phi\vee\psi)\Gamma(\eta\vee\zeta)$
\end{itemize}
The relation $\Gamma$ is thus obtained by propagating the constraint over $<\Theta>$.
It is obviously reflexive, symmetric and transitive; it is an equivalence relation.
An equivalence class for $\Gamma$ contains propositions which are identical in regards to the constraints.
\\[5pt]
It is now time to define the pre-Boolean algebra.
\paragraph{Pre-Boolean algebra.}
\begin{proposition}\label{dsmt2:dmb:prop3}
Let be given a free pre-Boolean algebra $<\Theta>$ and a set of propositional pairs \mbox{$\Gamma\subset<\Theta>\times<\Theta>$}\,.
Then, there is a magma denoted $<\Theta>_\Gamma$ and a morphism $\mu:<\Theta>\rightarrow<\Theta>_\Gamma$ such that:
$$
\left\{\begin{array}{l}
\mu\bigl(<\Theta>\bigr)=<\Theta>_\Gamma\;,
\vspace{3pt}\\
\forall\phi,\psi\in<\Theta>,\;\mu(\phi)=\mu(\psi)\iff\phi\Gamma\psi\;.
\end{array}\right.
$$
The magma $<\Theta>_\Gamma$ is called the pre-Boolean algebra generated by $\Theta$ and constrained by the constraints $\phi=\psi$ where $(\phi,\psi)\in\Gamma$\,.
\end{proposition}
\begin{description}
\item[Proof.]
For any $\phi\in<\Theta>$, define $\phi_\Gamma= \bigl\{\psi\in<\Theta>\; \big/\; \psi\Gamma\phi\bigr\}$; this set is called the class of $\phi$ for $\Gamma$\,.
\\
It is a well known fact, and the proof is immediate, that $\phi_\Gamma=\psi_\Gamma$ or $\phi_\Gamma\cap\psi_\Gamma=\emptyset$ for any $\phi,\psi\in<\Theta>$\,;
in particular, $\phi_\Gamma=\psi_\Gamma\iff\phi\Gamma\psi$\,.
\\
Now, assume $\eta_\Gamma=\phi_\Gamma$ and $\zeta_\Gamma=\psi_\Gamma$, that is $\eta\Gamma\phi$ and $\zeta\Gamma\psi$.\\
It comes $(\eta\wedge\zeta)\Gamma(\phi\wedge\psi)$ and $(\eta\vee\zeta)\Gamma(\phi\vee\psi)$.\\
As a consequence, $(\eta\wedge\zeta)_\Gamma=(\phi\wedge\psi)_\Gamma$ and $(\eta\vee\zeta)_\Gamma=(\phi\vee\psi)_\Gamma$.\\
At last:
$$
\Bigl(\eta_\Gamma=\phi_\Gamma\mbox{ and }
\zeta_\Gamma=\psi_\Gamma\Bigr)\Rightarrow
\Bigl((\eta\wedge\zeta)_\Gamma=(\phi\wedge\psi)_\Gamma\mbox{ and }
(\eta\vee\zeta)_\Gamma=(\phi\vee\psi)_\Gamma\Bigr)
$$
The proof is then concluded easily, by setting:
$$\left\{\begin{array}{@{}l@{}}
<\Theta>_\Gamma=\bigl\{\phi_\Gamma\;\big/\;\phi\in<\Theta>\bigr\}\;,
\vspace{3pt}\\
\forall\phi,\psi\in<\Theta>,\;
\phi_\Gamma\wedge\psi_\Gamma=(\phi\wedge\psi)_\Gamma
\mbox{ and }
\phi_\Gamma\vee\psi_\Gamma=(\phi\vee\psi)_\Gamma\;,
\vspace{3pt}\\
\forall\phi\in<\Theta>,\;\mu(\phi)=\phi_\Gamma\;.
\end{array}\right.
$$
\item[$\Box\Box\Box$]\rien
\end{description}
From now on, the element $\mu(\phi)$, where $\phi\in<\Theta>$, will be denoted $\phi$ as if $\phi$ were an element of $<\Theta>_\Gamma$\,.
In particular, $\mu(\phi)=\mu(\psi)$ will imply $\phi=\psi$ in $<\Theta>_\Gamma$ (but not in $<\Theta>$).
\begin{proposition}\label{dsmt2:dmb:prop4}
Let be given a free pre-Boolean algebra $<\Theta>$ and a set of propositional pairs \mbox{$\Gamma\subset<\Theta>\times<\Theta>$}\,.
Let $<\Theta>_\Gamma$ and $<\Theta>'_\Gamma$ be pre-Boolean algebras generated by $\Theta$ and constrained by the familly $\Gamma$\,.
Then $<\Theta>_\Gamma$ and $<\Theta>'_\Gamma$ are isomorph.
\end{proposition}
\begin{description}
\item[Proof.]
Let $\mu:<\Theta>\rightarrow<\Theta>_\Gamma$ and $\mu':<\Theta>\rightarrow<\Theta>'_\Gamma$ be as defined in proposition~\ref{dsmt2:dmb:prop3}.\\
For any $\phi\in<\Theta>$, define $\nu\bigl(\mu(\phi)\bigr)=\mu'(\phi)$\,.\\
Then, $\nu\bigl(\mu(\phi)\bigr)=\nu\bigl(\mu(\psi)\bigr)$ implies $\mu'(\phi)=\mu'(\psi)$\,.\\
By definition of $\mu'$, it is derived $\phi\Gamma\psi$ and then $\mu(\phi)=\mu(\psi)$\,.\\
Thus, $\nu$ is one-to-one.\\
By definition, it is also implied that $\nu$ is onto.
\item[$\Box\Box\Box$]\rien
\end{description}
This property thus says that there is a structural unicity of $<\Theta>_\Gamma$\,.
\begin{example}\label{DSMTb2:dmb:ex4}
Let us consider again the pre-Boolean algebra generated by \mbox{$\Theta=\{a,b,c\}$} and constrained by \mbox{$a\wedge b=a\wedge c$} and \mbox{$a\wedge c=b\wedge c$}.
In this case, the mapping $\mu:<\Theta>\rightarrow<\Theta>_\Gamma$ is defined by:
\begin{itemize}
\item $\mu\bigl(\{\bot\}\bigr)=\{\bot\}$, $\mu\Bigl(\bigl\{a,(b\wedge c)\vee a\bigr\}\Bigr)=\{a\}$, $\mu\Bigl(\bigl\{b,(c\wedge a)\vee b\bigr\}\Bigr)=\{b\}$, $\mu\Bigl(\bigl\{c,(a\wedge b)\vee c\bigr\}\Bigr)=\{c\}$, $\mu\bigl(\{a\vee b\vee c\}\bigr)=\{a\vee b\vee c\}$, $\mu\bigl(\{\top\}\bigr)=\{\top\}$\,,
\item $\mu\bigl(\{a\vee b\}\bigr)=\{a\vee b\}$, $\mu\bigl(\{b\vee c\}\bigr)=\{b\vee c\}$, $\mu\bigl(\{c\vee a\}\bigr)=\{c\vee a\}$\,,
\item $\mu\Bigl(\bigl\{a\wedge b\wedge c,a\wedge b,b\wedge c,c\wedge a,
(a\vee b)\wedge c, (b\vee c)\wedge a, (c\vee a)\wedge b,$
\\\rien\hspace{130pt}$(a\wedge b)\vee(b\wedge c)\vee(c\wedge a)
\bigr\}\Bigr)=\{a\wedge b\wedge c\}$\,.
\end{itemize}
\end{example}
\paragraph{Between sets and hyperpower sets.}
\begin{proposition}
The Boolean algebra $\bigl(\mathcal{P}(\Theta),\cap,\cup,\sim,\emptyset,\Theta\bigr)$, considered as a $(\wedge,\vee,\bot,\top)$-magma, is isomorph to the pre-Boolean algebra $<\Theta>_\Gamma$, where $\Gamma$ is defined by:
$$
\Gamma=\bigl\{(\theta\wedge\vartheta,\bot)\;\big/\;\theta,\vartheta\in\Theta\mbox{ and }\theta\ne\vartheta\bigr\}\cup\left\{\left(\bigvee_{\theta\in\Theta}\theta,\top\right)\right\}\;.
$$
\end{proposition}
\begin{description}
\item[Proof.]
Recall the notation $\varphi(\Sigma)=\bigvee_{\sigma\in\Sigma}\bigwedge_{\theta\in\sigma}\theta$ for any $\Sigma\subset\mathcal{P}(\Theta)$\,.\\
Define $\mu:<\Theta>\rightarrow\mathcal{P}(\Theta)$ by setting\footnote{It is defined $\bigcap_{\theta\in\emptyset}\theta=\Theta$\,.}
$\mu\bigl(\varphi(\Sigma)\bigr)=\bigcup_{\sigma\in\Sigma}\bigcap_{\theta\in\sigma}\{\theta\}$
for any $\Sigma\subset\mathcal{P}(\Theta)$\,.
\\
It is immediate that $\mu$ is a morphism.\\
Now, by definition of $\Gamma$, $\mu\bigl(\varphi(\Sigma)\bigr)=\mu\bigl(\varphi(\Lambda)\bigr)$ is equivalent to $\varphi(\Sigma)\Gamma\varphi(\Lambda)$\,.\\
The proof is then concluded by proposition~\ref{dsmt2:dmb:prop4}.
\item[$\Box\Box\Box$]\rien
\end{description}
Thus, sets, considered as Boolean algebra, and hyperpower sets are both extremal cases of the notion of pre-Boolean algebra.
But while hyperpower sets extend the structure of sets, hyperpower sets are more complexe in structure and size than sets.
A practical use of hyperpower sets becomes quickly impossible.
Pre-Boolean algebra however allows intermediate structures between sets and hyperpower sets.
\\[5pt]
A specific kind of pre-Boolean algebra will be particularly interesting when defining the DSmT.
Such pre-Boolean algebra will forbid any interaction between the trivial propositions $\bot,\top$ and the other propositions.
These algebra, called insulated pre-Boolean algebra, are characterized now.
\paragraph{Insulated pre-Boolean algebra.}
\label{f2k5:redefParagraph}
A pre-Boolean algebra $<\Theta>_\Gamma$ verifies the \emph{insulation} property if $\Gamma\subset\bigl(<\Theta>\setminus\{\bot,\top\})\bigr)\times\bigl(<\Theta>\setminus\{\bot,\top\})\bigr)$\,.
\begin{proposition}
Let $<\Theta>_\Gamma$ a pre-Boolean algebra verifying the insulation property.
Then holds for any $\phi,\psi\in<\Theta>_\Gamma$\,:
$$\left\{\begin{array}{l@{}}\displaystyle
\phi\wedge\psi=\bot\Rightarrow (\phi=\bot\mbox{ or }\psi=\bot)\;,
\vspace{3pt}\\\displaystyle
\phi\vee\psi=\top\Rightarrow (\phi=\top\mbox{ or }\psi=\top)
\;.
\end{array}\right.$$
In other words, all propositions are \emph{independent} with each other in a pre-Boolean algebra with insulation property.
\end{proposition}
The proof is immediate, since it is impossible to obtain $\phi\wedge\psi\Gamma\bot$ or $\phi\vee\psi\Gamma\top$ without involving $\bot$ or $\top$ in the constraints of $\Gamma$.
Examples~\ref{DSMTb2:dmb:ex2} and example~\ref{DSMTb2:dmb:ex3} verify the insulation property.
On the contrary, a non empty set does not.
%
\\[5pt]
\emph{Corollary and definition.}
Let $<\Theta>_\Gamma$ be a pre-Boolean algebra, verifying the insulation property.
Define $\ll\Theta\gg_\Gamma=<\Theta>_\Gamma\setminus\{\bot,\top\}$\,.
The operators $\wedge$ and $\vee$ restrict to $\ll\Theta\gg_\Gamma$\,,
and $\bigl(\ll\Theta\gg_\Gamma,\wedge,\vee\bigr)$
is an algebraic structure by itself, called \emph{insulated} pre-Boolean algebra.
This structure is also refered to as the insulated pre-Boolean algebra $\ll\Theta\gg_\Gamma$.
\begin{proposition}
Let $<\Theta>_\Gamma$ and $<\Theta>'_\Gamma$ be pre-Boolean algebras with insulation properties.
Assume that the insulated pre-Boolean algebra $\ll\Theta\gg_\Gamma$ and $\ll\Theta\gg'_\Gamma$ are $(\wedge,\vee)$-isomorph.
Then $<\Theta>_\Gamma$ and $<\Theta>'_\Gamma$ are isomorph.
\end{proposition}
Deduced from the insulation property.
\\\\
All ingredients are now gathered for the definition of Dezert Smarandache models.
\subsection{The \emph{free} Dezert Smarandache Theory}
\paragraph{Dezert Smarandache Model.}
Assume that $\Theta$ is a finite set.
A Dezert Smarandache model (DSmm) is a pair $(\Theta,m)$, where $\Theta$ is a set of propositions and the \emph{basic belief assignment} $m$ is a non negatively valued function defined over $<\Theta>$ such that:
$$
\sum_{\phi\in<\Theta>}m(\phi)=1
\quad\mbox{and}\quad
m(\bot)=0
\;.
$$
Moreover, it is generally assumed that $\sum_{\phi\in\ll\Theta\gg}m(\phi)=1$\,, which means that the propositions of $\Theta$ are exhaustive.
\paragraph{Belief Function.}
Assume that $\Theta$ is a finite set.
The belief function $\mathrm{Bel}$ related to a bba $m$ is defined by:
\begin{equation}
\label{f2k5:DSmTcont:Eq:1}
\forall\phi\in<\Theta>,\,\mathrm{Bel}(\phi)=\sum_{\psi\in<\Theta>:\psi\subset\phi}m(\psi)\;.
\end{equation}
The equation~(\ref{f2k5:DSmTcont:Eq:1}) is invertible:
$$
\forall\phi\in<\Theta>,\,m(\phi)=\mathrm{Bel}(\phi)-\sum_{\psi\in<\Theta>:\psi\subsetneq\phi}m(\psi)\;.
$$
\paragraph{Fusion rule.}
Assume that $\Theta$ is a finite set.
For a given universe $\Theta$\,, and two basic belief assignments $m_1$ and $m_2$, associated to independent sensors,
the fused basic belief assignment is $m_1\oplus m_2$\,, defined by:
\begin{equation}
\label{f2k5:DSmTcont:Eq:2}
m_1\oplus m_2(\phi)=\sum_{\psi_1,\psi_2\in<\Theta>:\psi_1\wedge\psi_2=\phi} m_1(\psi_1)m_2(\psi_2)\;.
\end{equation}
\paragraph{Remarks.}
It appears obviously that the previous definitions could be equivalently restricted to $\ll\Theta\gg$, owing to the insulation properties.
\\[3pt]
Considering the definition of the fusion rule and the insulation property $(\phi\ne\bot\mbox{ and }\psi\ne\bot)\Rightarrow(\phi\wedge\psi)\ne\bot$, it appears also that these definitions could be generalized to any algebra $<\Theta>_\Gamma$ with the insulation property.

\subsection{Extensions to any insulated pre-Boolean algebra}
\label{f2k5:subsectionExtentDSmT}
Let $\ll\Theta\gg_\Gamma$ be an insulated pre-Boolean algebra.
The definition of bba $m$, belief $\mathrm{Bel}$ and fusion $\oplus$ is thus kept unchanged (except the condition $m(\bot)=m(\top)=0$ which become useless).
\begin{itemize}
\item A \emph{basic belief assignment} $m$ is a non negatively valued function defined over $\ll\Theta\gg_\Gamma$ such that:
$$
\sum_{\phi\in\ll\Theta\gg_\Gamma}m(\phi)=1
\;.
$$
\item The belief function $\mathrm{Bel}$ related to a bba $m$ is defined by:
$$
\forall\phi\in\ll\Theta\gg_\Gamma,\,\mathrm{Bel}(\phi)=\sum_{\psi\in\ll\Theta\gg_\Gamma:\psi\subset\phi}m(\psi)\;.
$$
\item Being given two basic belief assignments $m_1$ and $m_2$,
the fused basic belief assignment $m_1\oplus m_2$ is defined by:
$$
m_1\oplus m_2(\phi)=\sum_{\psi_1,\psi_2\in\ll\Theta\gg_\Gamma:\psi_1\wedge\psi_2=\phi} m_1(\psi_1)m_2(\psi_2)\;.
$$
\end{itemize}
These extended definitions will be applied subsequently.
\section{Ordered DSm model}
\label{f2k5:sectionTowardOrderedModel}
From now on, we are working only with insulated pre-Boolean structures.
\\[5pt]
In order to reduce the complexity of the free DSm model, it is necessary to introduce logical constraints which will lower the size of the pre-Boolean algebra.
Such constraints may appear clearly in the hypotheses of the problem.
In this case, constraints come naturally and approximations may not be required.
However, when the model is too complex and there are no explicit constraints for reducing this complexity, it is necessary to approximate the model by introducing some new constraints.
Two rules should be applied then:
\begin{itemize}
\item Only weaken informations\footnote{Typically, a constraint like $\phi\wedge\psi\wedge\eta=\phi\wedge\psi$ will weaken the information, by erasing $\eta$ from $\phi\wedge\psi\wedge\eta$\,.}; do not produce information from nothing,
\item minimize the information weakening.
\end{itemize}
First point garantees that the approximation does not introduce false information.
But some significant informations (\emph{eg.} contradictions) are possibly missed.
This drawback should be avoided by second point.
\\[3pt]
In order to build a good approximation policy, some external knowledge, like distance or order relation among the propositions could be used.
Behind these relations will be assumed some kind of distance between the informations:
\emph{more are the informations distant, more are their conjunctive combination valuable.}
\subsection{Ordered atomic propositions}
Let $(\Theta,\le)$ be an ordered set of atomic propositions.
This order relation is assumed to describe the relative distance between the information.
For example, the relation $\phi\le\psi\le\eta$ implies that $\phi$ and $\psi$ are closer informations than $\phi$ and $\eta$\,.
Thus, the information contained in $\phi\wedge\eta$ is stronger than the information contained in $\phi\wedge\psi$\,.
Of course, this comparison does not matter when all the information is kept, but when approximations are necessary, it will be useful to be able to choose the best information.
\paragraph{Sketchy example.}
Assume that 3 independent sensors are giving $3$ measures about a continuous parameter, that is $x$, $y$ and $z$.
The parameters $x,y,z$ are assumed to be real values, not of the \emph{set} $\Rset$ but of its pre-Boolean extension (theoretical issues will be clarified later\footnote{In particular, as we are working in a pre-Boolean algebra, $x\wedge y$ makes sense and it is possible that $x\wedge y\ne\bot$ even when $x\ne y$\,.}%
).
The fused information could be formalized by the proposition $x\wedge y\wedge z$ (in a DSmT viewpoint).
What happen if we want to reduce the information by removing a proposition.
Do we keep $x\wedge y$\,, $y\wedge z$ or $x\wedge z$\,?
This is of course an information weakening.
But it is possible that one information is better than an other.
At this stage, the order between the values $x,y,z$ will be involved.
Assume for example that $x\le y<z$\,.
It is clear that the proposition $x\wedge z$ indicates a greater contradiction than $x\wedge y$ or $y\wedge z$\,.
Thus, the proposition $x\wedge z$ is the one which should be kept!
The discarding constraint $x\le y\le z \Rightarrow x\wedge y\wedge z=x\wedge z$ is implied then.
\subsection{Associated pre-Boolean algebra and complexity.}
In regard to the previous example, the insulated pre-Boolean algebra associated to the ordered propositions $(\Theta,\le)$ is $\ll\Theta\gg_{\Gamma}$\,, where $\Gamma$ is defined by:
$$
\Gamma=\bigl\{(\phi\wedge\psi\wedge\eta,\phi\wedge\eta)\big/\phi,\psi,\eta\in\Theta\mbox{ and }\phi\le\psi\le\eta\bigr\}\;.
$$
The following property give an approximative bound of the size of $\ll\Theta\gg_{\Gamma}$ in the case of a total order.
\begin{proposition}
Assume that $(\Theta,\le)$ is totally ordered.
Then, $\ll\Theta\gg_{\Gamma}$ is a substructure of the set $\Theta^{2}$\,.
\end{proposition}
\begin{description}
\item[proof.] 
Since the order is total, first notice that the added constraints are:
$$
\forall \phi,\psi,\eta\in\Theta\,,\;\phi\wedge\psi\wedge\eta=\min\{\phi,\psi,\eta\}\wedge\max\{\phi,\psi,\eta\}\;.
$$
Now, for any $\phi\in\Theta$\,, define $\breve{\phi}$ by:
$$
\underline{\breve{\phi}=\bigl\{(\varphi_1,\varphi_2)\in\Theta^2\,\big/\,\varphi_1\le\phi\le\varphi_2\bigr\}
}
$$
It is noteworthy that:
$$
\breve{\phi}\cap\breve{\psi}=\bigl\{(\varphi_1,\varphi_2)\in\Theta^2\,\big/\,\varphi_1\le\min\{\phi,\psi\}\mbox{ and }\max\{\phi,\psi\}\le\varphi_2\bigr\}
$$
and
$$
\breve{\phi}\cap\breve{\psi}\cap\breve{\eta}=\bigl\{(\varphi_1,\varphi_2)\in\Theta^2\,\big/\,\varphi_1\le\min\{\phi,\psi,\eta\}\mbox{ and }\max\{\phi,\psi,\eta\}\le\varphi_2\bigr\}\;.
$$
By defining $m=\min\{\phi,\psi,\eta\}$ and $M=\max\{\phi,\psi,\eta\}$\,, it is deduced:
\begin{equation}\label{DSmTbk2:inPrf:1}
\breve{\phi}\cap\breve{\psi}\cap\breve{\eta}=\breve{m}\cap\breve{M}\;.
\end{equation}
Figure~\ref{DSmTb2:dmb:fir:ord:1} illustrates the construction of $\breve{\phi}$, $\breve{\phi}\cap\breve{\psi}$ and property~(\ref{DSmTbk2:inPrf:1}).
\begin{figure}[h!]
\centering
\scalebox{1}{\begin{tabular}{@{}c@{ \ }c@{ \ }c@{}}
\xy
(0,0)*+{\rien},
(0,-24)="fromX",
(32,-24)*+{\varphi_1}="toX",
(8,-32)="fromY",
(8,0)*+{\varphi_2}="toY",
(0,-32)="fromXY",
(32,0)="toXY",
(30,-8)*+{_{\varphi_2=\varphi_1}},
(16,-24)="phi",
(16,-28)*+{\phi},
(16,-16)="xyphi",
(0,-16)="xphi",
(16,0)="yphi",
(2,-12)*+{\breve\phi},
\ar @{--}"phi";"xyphi",
\ar @{=}"xyphi";"xphi",
\ar @{=}"xyphi";"yphi",
\ar @{->}"fromX";"toX",
\ar @{->}"fromY";"toY",
\ar @{--}"toXY";"fromXY",
\POS "fromXY",
%
(32,-32)*+{\rien}
\endxy
&
\xy
(0,0)*+{\rien},
(0,-24)="fromX",
(32,-24)="toX",
(8,-32)="fromY",
(8,0)="toY",
(0,-32)="fromXY",
(32,0)="toXY",
(12,-24)="phi",
(12,-28)*+{\phi},
(12,-20)="xyphi",
(0,-20)="xphi",
(12,0)="yphi",
(20,-24)="psi",
(20,-28)*+{\psi},
(20,-12)="xypsi",
(0,-12)="xpsi",
(20,0)="ypsi",
(2,-8)*+{\breve\phi\cap\breve\psi},
(12,-12)="xyphipsi",
\ar @{--}"phi";"xyphi",
\ar @{:}"xyphi";"xphi",
\ar @{:}"xyphi";"xyphipsi",
\ar @{=}"xyphipsi";"yphi",
\ar @{--}"psi";"xypsi",
\ar @{:}"xypsi";"xyphipsi",
\ar @{=}"xyphipsi";"xpsi",
\ar @{:}"xypsi";"ypsi",
\ar @{->}"fromX";"toX",
\ar @{->}"fromY";"toY",
\ar @{--}"toXY";"fromXY",
\POS "fromXY",
%
(32,-32)*+{\rien}
\endxy
&
\xy
(0,0)*+{\rien},
(0,-32)="fromXY",
(32,0)="toXY",
(15,-26)*+{\breve\phi=\breve m},
(8,-24)="xyphi",
(0,-24)="xphi",
(8,0)="yphi",
(18,-18)*+{\breve\eta},
(16,-16)="xyeta",
(0,-16)="xeta",
(16,0)="yeta",
(31,-10)*+{\breve\psi=\breve M},
(24,-8)="xypsi",
(0,-8)="xpsi",
(24,0)="ypsi",
(0,-4)*+{\breve m\cap\breve M},
(8,-8)="xyphipsi",
\ar @{:}"xyphi";"xphi",
\ar @{:}"xyphi";"xyphipsi",
\ar @{=}"xyphipsi";"yphi",
\ar @{:}"xyeta";"xeta",
\ar @{:}"xyeta";"yeta",
\ar @{:}"xypsi";"xyphipsi",
\ar @{=}"xyphipsi";"xpsi",
\ar @{:}"xypsi";"ypsi",
\ar @{--}"toXY";"fromXY",
\POS "fromXY",
%
(32,-32)*+{\rien}
\endxy
\\
$\breve\phi$
&
$\breve\phi\cap\breve\psi$
&
$\breve\phi\cap\breve\psi\cap\breve\eta=\breve m\cap\breve M$
\end{tabular}}
\caption{Construction of $\breve\phi$}
\label{DSmTb2:dmb:fir:ord:1}
\end{figure}
\\[5pt]
Let $\mathcal{A}\subset\mathcal{P}(\Theta^2)$ be generated by $\breve{\phi}|_{\phi\in\Theta}$ with $\cap$ and $\cup$\,, \emph{ie.}:
$$
\mathcal{A}=\bigcup_{n\ge0}\left\{\left.\bigcup_{k=1}^n\bigl(\breve{\phi}_k\cap\breve{\psi}_k\bigr)\right/
\forall k\,,\;\breve{\phi}_k,\breve{\psi}_k\in\Theta\right\}\;.
$$
A consequence of~(\ref{DSmTbk2:inPrf:1}) is that $\mathcal{A}$ is an insulated pre-Boolean algebra which satisfies the constraints of $\Gamma$.
Then, the mapping:
$$\fbox{$\displaystyle
\smallsmile\;:\left\{\begin{array}{@{}l@{}}\displaystyle
\ll\Theta\gg_\Gamma\;\;\longrightarrow\;\mathcal{A}
\\\displaystyle
\bigvee_{k=1}^n\bigwedge_{l=1}^{n_k}\phi_{k,l}
\longmapsto
\bigcup_{k=1}^n\bigcap_{l=1}^{n_k}{\breve{\phi}}_{k,l}\quad,\quad\mbox{where }\phi_{k,l}\in\Theta
\end{array}\right.
$}$$
is an onto morphism of pre-Boolean algebra.\\
Now, let us prove that $\smallsmile$ is a one-to-one morphism.
\begin{lemma}
Assume:
$$
\bigcup_{k=1}^n\bigl(\breve{\phi}_k^1\cap\breve{\phi}_k^2\bigr)
\subset
\bigcup_{l=1}^m\bigl(\breve{\psi}_{l}^1\cap\breve{\psi}_{l}^2\bigr)
\quad,\quad\mbox{where }\phi_{k}^j\;,\;\psi_{l}^j\in\Theta\;.
$$
Then:
$$
\forall k\,,\,\exists l\,,\;\min\{\phi_k^1,\phi_k^2\}\le\min\{\psi_l^1,\psi_l^2\}
\ \mbox{and}\ 
\max\{\phi_k^1,\phi_k^2\}\ge\max\{\psi_l^1,\psi_l^2\}
$$
and
$$
\forall k\,,\,\exists l\,,\;\breve{\phi}_k^1\cap\breve{\phi}_k^2\subset\breve{\psi}_{l}^1\cap\breve{\psi}_{l}^2\;.
$$
\end{lemma}
\begin{description}
\item[Proof of lemma.]
Let $k\in[\![1,n]\!]$\,.\\
Define $m=\min\{\phi_k^1,\phi_k^2\}$ and $M=\max\{\phi_k^1,\phi_k^2\}$\,.\\
Then holds $(m,M)\in\breve{\phi}_k^1\cap\breve{\phi}_k^2$\,, implying $(m,M)\in\bigcup_{l=1}^m\bigl(\breve{\psi}_{l}^1\cap\breve{\psi}_{l}^2\bigr)$\,.\\
Let $l$ be such that $(m,M)\in\breve{\psi}_{l}^1\cap\breve{\psi}_{l}^2$\,.\\ 
Then $m\le\min\{\psi_l^1,\psi_l^2\}$ and $M\ge\max\{\psi_l^1,\psi_l^2\}$\,.\\
Followingly, $\breve{\phi}_k^1\cap\breve{\phi}_k^2\subset\breve{\psi}_{l}^1\cap\breve{\psi}_{l}^2$\,.
\item[$\Box\Box$]\rien
\end{description}
From $\min\{\phi_k^1,\phi_k^2\}\le\min\{\psi_l^1,\psi_l^2\}$
and 
$\max\{\phi_k^1,\phi_k^2\}\ge\max\{\psi_l^1,\psi_l^2\}$ is also deduced $(\phi_k^1\wedge\phi_k^2)\wedge(\psi_l^1\wedge\psi_l^2)=\phi_k^1\wedge\phi_k^2$ (definition of $\Gamma$)\,.\\
This property just means $\phi_k^1\wedge\phi_k^2\subset\psi_l^1\wedge\psi_l^2$\,.
It is lastly deduced:
\begin{lemma}
Assume:
$$
\bigcup_{k=1}^n\bigl(\breve{\phi}_k^1\cap\breve{\phi}_k^2\bigr)
\subset
\bigcup_{l=1}^m\bigl(\breve{\psi}_{l}^1\cap\breve{\psi}_{l}^2\bigr)
\quad,\quad\mbox{where }\phi_{k}^j\;,\;\psi_{l}^j\in\Theta\;.
$$
Then:
$$
\bigvee_{k=1}^n\bigl(\phi_k^1\wedge\phi_k^2\bigr)
\subset
\bigvee_{l=1}^m\bigl(\psi_{l}^1\wedge\psi_{l}^2\bigr)
\;.
$$
\end{lemma}
From this lemma, it is deduced that $\smallsmile$ is one to one.\\
At last \underline{$\smallsmile$ is an isomorphism of pre-Boolean algebra}, and $\ll\Theta\gg_\Gamma$ is a substructure of $\Theta^2$\,.
\item[$\Box\Box\Box$]\rien
\end{description}
\subsection{General properties of the model}
\begin{figure}[ht!]
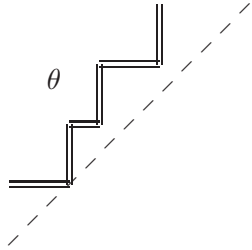

\centering
\scalebox{1}{\begin{tabular}{@{}c@{}}
\xy
(0,0)*+{\rien},
(0,-32)="fromXY",
(32,0)="toXY",
(0,-24)="a",
(8,-24)="b",
(8,-16)="c",
(12,-16)="d",
(12,-8)="e",
(20,-8)="f",
(20,0)="g",
(6,-10)*+{\theta},
\ar @{=}"a";"b",
\ar @{=}"b";"c",
\ar @{=}"c";"d",
\ar @{=}"d";"e",
\ar @{=}"e";"f",
\ar @{=}"f";"g",
\ar @{--}"toXY";"fromXY",
\POS "fromXY",
%
(32,-32)*+{\rien}
\endxy
\end{tabular}}
\caption{Example of increasing subset $\theta\in\mathcal{U}$}
\label{DSmTb2:dmb:fir:ord:2}
\end{figure}

In the next section, the previous construction will be extended to the continuous case, \emph{ie.} $(\Rset,\le)$\,.
However, a strict logical manipulation of the propositions is not sufficient and instead a measurable generalization of the model will be used.
It has been seen that a proposition of $\ll\Theta\gg_\Gamma$ could be described as a subset of $\Theta^2$\,.
In this subsection, the proposition model will be characterized precisely.
This characterization will be used and extended in the next section to the continuous case.
\begin{proposition}
Let $\phi\in\ll\Theta\gg_\Gamma$\,.\\
Then $\smallsmile\!(\phi)\subset\mathcal{T}$\,, where $\mathcal{T}=\bigl\{(\phi,\psi)\in\Theta^2\big/\phi\le\psi\bigr\}$\,.
\end{proposition}
\begin{description}
\item[Proof.] Obvious, since $\forall\phi\in\Theta\,,\;\breve{\phi}\subset\mathcal{T}$\,.
\item[$\Box\Box\Box$]\rien
\end{description}
\begin{definition}
A subset $\theta\subset\Theta^2$ is \emph{increasing} if and only if:
$$
\forall\,(\phi,\psi)\in\theta\,,\;\forall\eta\le\phi\,,\;\forall\zeta\ge\psi\,,\;(\eta,\zeta)\in\theta\;.
$$
\end{definition}
Let $\mathcal{U}=\bigl\{\theta\subset\mathcal{T}\big/\theta\mbox{ is increasing}\mbox{ and }\theta\ne\emptyset\bigr\}$ be the set of increasing non-empty subsets of $\mathcal{T}$\,.
Notice that the intersection or the union of increasing non-empty subsets are increasing non-empty subsets, so that $(\mathcal{U},\cap,\cup)$ is an insulated pre-Boolean algebra.
\begin{proposition}
For any choice of $\Theta$\,, $\bigl\{\smallsmile\!(\phi)\big/\phi\in\ll\Theta\gg_{\Gamma}\bigr\}\subset\mathcal{U}\;$.
\\
When $\Theta$ is finite, $\mathcal{U}=\bigl\{\smallsmile\!(\phi)\big/\phi\in\ll\Theta\gg_{\Gamma}\bigr\}\;$.
\end{proposition}
\begin{description}
\item[Proof of $\supset$\,.]
Obvious, since $\breve{\phi}$ is inceasing for any $\phi\in\Theta$\,.
\item[Proof of $\subset$\,.]
Let $\theta\in\mathcal{U}$ and let $(a,b)\in\theta$\,.\\
Since $\breve{a}\cap\breve{b}=\bigl\{(\alpha,\beta)\in\Theta^2\big/\alpha\le a\mbox{ and }\beta\ge b\bigr\}$ and $\theta$ is increasing, it follows $\breve{a}\cap\breve{b}\subset\theta$\,.\\
At last, $\theta=\bigcup_{(a,b)\in\theta}\breve{a}\cap\breve{b}=\smallsmile\!\left(\bigvee_{(a,b)\in\theta}a\wedge b\right)$\,.\\
Notice that $\bigvee_{(a,b)\in\theta}a\wedge b$ is actually defined, since $\theta$ is finite when $\Theta$ is finite.
\item[$\Box\Box\Box$]\rien
\end{description}
Figure~\ref{DSmTb2:dmb:fir:ord:2} gives an example of increasing subsets, element of $\mathcal{U}$.\\[5pt]
When infinite $\vee$-ing are allowed, notice that $\mathcal{U}$ may be considered as a model for $\ll\Theta\gg_{\Gamma}$ even if $\Theta$ is infinite.
In the next section, the \emph{continuous} pre-Boolean algebra related to $(\Rset,\le)$ will be modelled by the \emph{measurable increasing subsets} of \mbox{$\bigl\{(x,y)\in\Rset^2\big/x\le y\bigr\}$}\,.
\section{Continuous DSm model}
\label{f2k5:sectionTowardContinuousModel}
In this section, the case $\Theta=\Rset$ is considered.
\\[5pt]
Typically, in a continuous model, it will be necessary to manipulate any measurable proposition, and for example intervals.
It comes out that most intervals could not be obtained by a finite logical combinaison of the atomic propositions, but rather by infinite combinations.
For example, considering the set formalism, it is obtained $[a,b]=\bigcup_{x\in[a,b]}\{x\}$\,, which suggests the definition of the infinite disjunction ``$\bigvee_{x\in[a,b]}x$''.
It is known that infinite disjunctions are difficult to handle in a logic.
It is better to manipulate the models directly.
The pre-Boolean algebra to be constructed should verify the property $x\le y\le z \Rightarrow x\wedge y\wedge z=x\wedge z$\,.
As discussed previously and since infinitary disjunctions are allowed, a model for such algebra are the measurable increasing subsets.
\subsection{Measurable increasing subsets}
\begin{figure}[ht!]
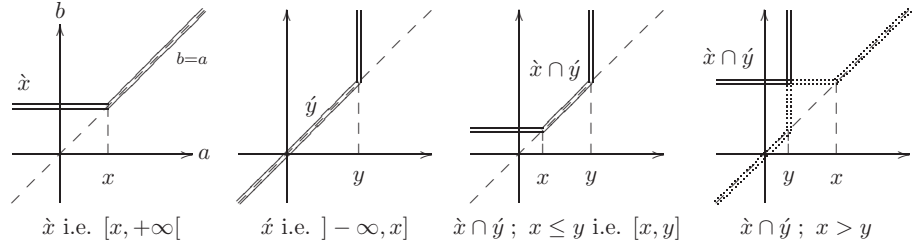

\centering
\scalebox{.8}{\begin{tabular}{@{}c@{ \ }c@{ \ }c@{ \ }c@{}}
\xy
(0,0)*+{\rien},
(0,-24)="fromX",
(32,-24)*+{a}="toX",
(8,-32)="fromY",
(8,0)*+{b}="toY",
(0,-32)="fromXY",
(32,0)="toXY",
(30,-8)*+{_{b=a}},
(16,-24)="phi",
(16,-28)*+{x},
(16,-16)="xyphi",
(0,-16)="xphi",
(2,-12)*+{\grave x},
\ar @{--}"phi";"xyphi",
\ar @{=}"xyphi";"xphi",
\ar @{=}"xyphi";"toXY",
\ar @{->}"fromX";"toX",
\ar @{->}"fromY";"toY",
\ar @{--}"toXY";"fromXY",
\POS "fromXY",
%
(32,-32)*+{\rien}
\endxy
&
\xy
(0,0)*+{\rien},
(0,-24)="fromX",
(32,-24)="toX",
(8,-32)="fromY",
(8,0)="toY",
(0,-32)="fromXY",
(32,0)="toXY",
(20,-24)="psi",
(20,-28)*+{y},
(20,-12)="xypsi",
(20,0)="ypsi",
(12,-16)*+{\acute y},
%
\ar @{--}"psi";"xypsi",
\ar @{=}"xypsi";"ypsi",
\ar @{=}"xypsi";"fromXY",
%
%
\ar @{->}"fromX";"toX",
\ar @{->}"fromY";"toY",
\ar @{--}"toXY";"fromXY",
\POS "fromXY",
%
(32,-32)*+{\rien}
\endxy
&
\xy
(0,0)*+{\rien},
(0,-24)="fromX",
(32,-24)="toX",
(8,-32)="fromY",
(8,0)="toY",
(0,-32)="fromXY",
(32,0)="toXY",
(12,-24)="phi",
(12,-28)*+{x},
(12,-20)="xyphi",
(0,-20)="xphi",
%
(20,-24)="psi",
(20,-28)*+{y},
(20,-12)="xypsi",
(20,0)="ypsi",
(14,-10)*+{\grave x\cap\acute y},
\ar @{--}"phi";"xyphi",
\ar @{--}"psi";"xypsi",
\ar @{=}"xyphi";"xphi",
\ar @{=}"xyphi";"xypsi",
\ar @{=}"xypsi";"ypsi",
\ar @{->}"fromX";"toX",
\ar @{->}"fromY";"toY",
\ar @{--}"toXY";"fromXY",
\POS "fromXY",
%
(32,-32)*+{\rien}
\endxy
&
\xy
(0,0)*+{\rien},
(0,-24)="fromX",
(32,-24)="toX",
(8,-32)="fromY",
(8,0)="toY",
(0,-32)="fromXY",
(32,0)="toXY",
(12,-24)="phi",
(12,-28)*+{y},
(12,-20)="xyphi",
(12,0)="yphi",
(20,-24)="psi",
(20,-28)*+{x},
(20,-12)="xypsi",
(0,-12)="xpsi",
(2,-8)*+{\grave x\cap\acute y},
(12,-12)="xyphipsi",
\ar @{--}"phi";"xyphi",
\ar @{:}"xyphi";"fromXY",
\ar @{:}"xyphi";"xyphipsi",
\ar @{=}"xyphipsi";"yphi",
\ar @{--}"psi";"xypsi",
\ar @{:}"xypsi";"xyphipsi",
\ar @{=}"xyphipsi";"xpsi",
\ar @{:}"xypsi";"toXY",
\ar @{->}"fromX";"toX",
\ar @{->}"fromY";"toY",
\ar @{--}"toXY";"fromXY",
\POS "fromXY",
%
(32,-32)*+{\rien}
\endxy
\\
$\grave x\mbox{ i.e. }[x,+\infty[$
&
$\acute x\mbox{ i.e. }]-\infty,x]$
&
$\grave x\cap\acute y\;;\ x\le y\mbox{ i.e. }[x,y]$
&
$\grave x\cap\acute y\;;\ x>y$
\end{tabular}}
\caption{Interval models}
\label{DSmTb2:dmb:fir:ord:3}
\end{figure}

A measurable subset $A\subset \Rset^2$ is a measurable increasing subset if:
$$\left\{\begin{array}{@{}l@{}}\displaystyle
\forall\,(x,y)\in A\,,\; x\le y\;,
\\\displaystyle
\forall\,(x,y)\in A\,,\;\forall a\le x\,,\;\forall b\ge y\,,\;(a,b)\in A\;.
\end{array}\right.$$
The set of measurable increasing subsets is denoted $\mathcal{U}$.
\paragraph{Example.}
Let $f:\Rset\rightarrow\Rset$ be a non decreasing measurable mapping such that $f(x)\ge x$ for any $x\in\Rset$.
The set $\bigl\{(x,y)\in\Rset^2\big/f(x)\le y\bigr\}$ is a measurable increasing subset.
\paragraph{``Points''.}
For any $x\in\Rset$, the measurable increasing subset $\breve{x}$ is defined by:
$$
\breve{x}=\bigl\{(a,b)\in\Rset^2\ \big/\ a\le x\le b\bigr\}\;.
$$
The set $\breve{x}$ is of course a model for the point $x\in\Rset$ within the pre-Boolean algebra (refer to section~\ref{f2k5:sectionTowardOrderedModel}).
\paragraph{Generalized intervals.}
A particular class of increasing subsets, the generalized intervals, will be useful in the sequel.
\\[5pt]
For any $x\in\Rset$, the measurable sets $\grave{x}$ and $\acute{x}$ are defined by:
$$
\left\{\begin{array}{@{\,}l@{}}\vspace{5pt}\displaystyle
\grave{x}=\bigl\{(a,b)\in\Rset^2\ \big/\ a\le b\mbox{ and }x\le b\bigr\}\;,
\\\displaystyle
\acute{x}=\bigl\{(a,b)\in\Rset^2\ \big/\ a\le b\mbox{ and }a\le x\bigr\}\;.
\end{array}\right.
$$
\\[5pt]
The following properties are derived:
$$
\breve{x}=\grave{x}\cap\acute{x}\;,
\ 
\grave{x}=\bigcup_{z\in[x,+\infty[}\breve{z}
\quad\mbox{and}\quad
\acute{x}=\bigcup_{z\in]-\infty,x]}\breve{z}
$$
Moreover, for any $x,y$ such that $x\le y$, it comes:
$$
\grave{x}\cap\acute{y}=\bigcup_{z\in [x,y]}\breve{z}\;.
$$
As a conclusion, the set $\grave{x}$, $\acute{x}$ and $\grave{x}\cap\acute{y}$ (with $x\le y$) are the respective models for the intervals $[x,+\infty[$\,, $]-\infty,x]$ and $[x,y]$ within the pre-Boolean algebra.
Naturally, the quotation marks $\grave{}$ (opening) and $\acute{}$ (closing) are used respectively for opening and closing the intervals.
Figure~\ref{DSmTb2:dmb:fir:ord:3} illustrates various cases of interval models.
\\[5pt]
At last, the set $\grave{x}\cap\acute{y}$, where $x,y\in\Rset$ are not constrained, constitutes a generalized definition of the notion of interval.
In the case $x\le y$, it works like ``classical'' interval, but in the case $x>y$, it is obtained a new class of intervals with negative width (last case in figure~\ref{DSmTb2:dmb:fir:ord:3}).
Whatever, $\grave{x}\cap\acute{y}$ comes with a non empty inner, and may have a non zero measure.
\\[3pt]
The width $\delta=\frac{y-x}2$ of the interval $\grave{x}\cap\acute{y}$ could be considered as a measure of contradiction associated with this proposition, while its center $\mu=\frac{x+y}2$ should be considered as its median value.
The interpretation of the measure of contradiction is left to the human.
Typically, a possible interpretation could be:
\begin{itemize}
\item $\delta<0$ means contradictory informations,
\item $\delta=0$ means exact informations,
\item $\delta>0$ means imprecise informations.
\end{itemize}
It is also noteworthy that the set of generalized intervals $$\underline{\mathcal{I}=\left\{\grave{x}\cap\acute{y}/x,y\in\Rset\right\}}$$
is left unchanged by the operator $\cap$\,, as seen in the following proposition~\ref{f2k5:lemmaStability}\,:
\begin{proposition}[Stability]\label{f2k5:lemmaStability}
\emph{Let $x_1,x_2,y_1,y_2\in\Rset$\,.\\
Define \mbox{$x=\max\{x_1,x_2\}$} and \mbox{$y=\min\{y_1,y_2\}$}\,.\\
Then
$\bigl(\grave{x}_1\cap\acute{y}_1\bigr)
\cap\bigl(\grave{x}_2\cap\acute{y}_2\bigr)=\grave{x}\cap\acute{y}$
\,.}
\end{proposition}
\emph{Proof is obvious.}
\\[5pt]
This last property make possible the definition of basic belief assignment over generalized intervals only.
This assumption is clearly necessary in order to reduce the complexity of the evidence modelling.
Behind this assumption is the idea that a continuous measure is described by an  imprecision/contradiction around the measured value.
Such hypothesis has been made by Smets and Ristic\cite{risticsmets}.
From now on, all the defined bba will be zeroed outside $\mathcal{I}$.
Now, since $\mathcal{I}$ is invariant by $\cap$\,, it is implied that all the bba which will be manipulated, from sensors or after fusion, will be zeroed outside $\mathcal{I}$.
This makes the basic belief assignments equivalent to a density over the 2-dimension space $\Rset^2$\,.
\subsection{Definition and manipulation of the belief}
The definitions of bba, belief and fusion result directly from section~\ref{f2k5:sectionIntroDSmT}, but of course the bba becomes density and the summations are replaced by integrations.
\paragraph{Basic Belief Assignment.}
As discussed previously, it is hypothesized that the measures are characterized by a precision interval around the measured values.
In addition, there is an uncertainty about the measure which is translated into a basic belief assignment over the precision intervals.
\\[5pt]
According to these hypotheses, a bba will be a non negatively valued function $m$ defined over $\mathcal{U}$\,, zeroed outside $\mathcal{I}$ (set of generalized intervals), and such that:
$$
\int_{x,y\in\Rset}m\bigl(\grave{x}\cap\acute{y}\bigr)dxdy=1\;.
$$
\paragraph{Belief function.}
The function of belief, $\mathrm{Bel}$, is defined for any measurable proposition $\phi\in\mathcal{U}$ by:
$$
\mathrm{Bel}\,(\phi)=\int_{\grave{x}\cap\acute{y}\subset\phi} m\bigl(\grave{x}\cap\acute{y}\bigr)dxdy\;.
$$
In particular, for a generalized interval $\grave{x}\cap\acute{y}$\,:
$$
\mathrm{Bel}\,\bigl(\grave{x}\cap\acute{y}\bigr)=\int_{u=x}^{+\infty}\int_{v=-\infty}^{y} m\bigl(\grave{u}\cap\acute{v}\bigr)dudv\;.
$$
\paragraph{Fusion rule.}
Being given two basic belief assignments $m_1$ and $m_2$,
the fused basic belief assignment $m_1\oplus m_2$ is defined by the curviline integral:
$$
m_1\oplus m_2\bigl(\grave{x}\cap\acute{y}\bigr)=\int_{\mathcal{C}=\{(\phi,\psi)/\phi\cap\psi=\grave{x}\cap\acute{y}\}}m_1(\phi)m_2(\psi)\,d\mathcal{C}\;.
$$
Now, from hypothesis it is assumed that $m_i$ is positive only for intervals of the form $\grave{x}_i\cap\acute{y}_i$.
Proposition~\ref{f2k5:lemmaStability} implies:
$$
\grave{x}_1\cap\acute{y}_1\cap\grave{x}_2\cap\acute{y}_2=\grave{x}\cap\acute{y}
\mbox{ where }
\left\{\begin{array}{@{\,}l}\vspace{3pt}\displaystyle
x=\max\{x_1,x_2\}\;,
\\\displaystyle
y=\min\{y_1,y_2\}\;.
\end{array}\right.
$$
It is then deduced:
$$\begin{array}{@{\,}l}\vspace{3pt}\displaystyle
m_1\oplus m_2\bigl(\grave{x}\cap\acute{y}\bigr)=
\int_{x_2=-\infty}^{x}\int_{y_2=y}^{+\infty}
m_1\bigl(\grave{x}\cap\acute{y}\bigr)m_2\bigl(\grave{x}_2\cap\acute{y}_2\bigr)dx_2dy_2
\\\vspace{3pt}\displaystyle\hspace{30pt}
+\int_{x_1=-\infty}^{x}\int_{y_1=y}^{+\infty}
m_1\bigl(\grave{x}_1\cap\acute{y}_1\bigr)m_2\bigl(\grave{x}\cap\acute{y}\bigr)dx_1dy_1
\\\vspace{3pt}\displaystyle\hspace{30pt}
+\int_{x_1=-\infty}^{x}\int_{y_2=y}^{+\infty}
m_1\bigl(\grave{x}_1\cap\acute{y}\bigr)m_2\bigl(\grave{x}\cap\acute{y}_2\bigr)dx_1dy_2
\\\displaystyle\hspace{30pt}
+\int_{x_2=-\infty}^{x}\int_{y_1=y}^{+\infty}
m_1\bigl(\grave{x}\cap\acute{y}_1\bigr)m_2\bigl(\grave{x}_2\cap\acute{y}\bigr)dx_2dy_1
\;.
\end{array}$$
In particular, it is now justified that a bba, from sensors or fused, will always be zeroed outside $\mathcal{I}$\,.
\section{Implementation of the continuous model}
\label{f2k5:sectionImplementation}
\begin{figure}
\centering
\scalebox{.75}{\input{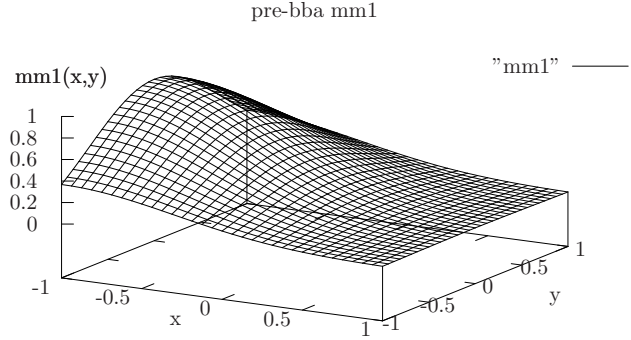}}
\caption{Non normalized bba $mm_1$}
\label{f2k5:fig:mm1}
\end{figure}

\paragraph{Setting.}
In this implementation, the study has been restricted to the interval $[-1,1]$ instead of $\Rset$.
The previous results still hold by trunctating over $[-1,1]$\,.
In particular, any bba $m$ is zeroed outside $\mathcal{I}_{-1}^{1}=\left\{\grave{x}\cap\acute{y}/x,y\in[-1,1]\right\}$ and its related belief function is defined by:
$$
\mathrm{Bel}\,\bigl(\grave{x}\cap\acute{y}\bigr)=\int_{u=x}^{1}\int_{v=-1}^{y} m\bigl(\grave{u}\cap\acute{v}\bigr)dudv\;,
$$
for any generalized interval of $\mathcal{I}_{-1}^{1}$\,.
The bba resulting of the fusion of two bbas $m_1$ and $m_2$ is defined by:
$$\begin{array}{@{\,}l}\vspace{3pt}\displaystyle
m_1\oplus m_2\bigl(\grave{x}\cap\acute{y}\bigr)=
\int_{x_2=-1}^{x}\int_{y_2=y}^{1}
m_1\bigl(\grave{x}\cap\acute{y}\bigr)m_2\bigl(\grave{x}_2\cap\acute{y}_2\bigr)dx_2dy_2
\\\vspace{3pt}\displaystyle\hspace{30pt}
+\int_{x_1=-1}^{x}\int_{y_1=y}^{1}m_1\bigl(\grave{x}_1\cap\acute{y}_1\bigr)m_2\bigl(\grave{x}\cap\acute{y}\bigr)dx_1dy_1
\\\vspace{3pt}\displaystyle\hspace{30pt}
+\int_{x_1=-1}^{x}\int_{y_2=y}^{1}m_1\bigl(\grave{x}_1\cap\acute{y}\bigr)m_2\bigl(\grave{x}\cap\acute{y}_2\bigr)dx_1dy_2
\\\displaystyle\hspace{30pt}
+\int_{x_2=-1}^{x}\int_{y_1=y}^{1}m_1\bigl(\grave{x}\cap\acute{y}_1\bigr)m_2\bigl(\grave{x}_2\cap\acute{y}\bigr)dx_2dy_1
\;.
\end{array}$$
\paragraph{Method.}
A theorical computation of these integrals seems uneasy.
An approximation of the densities and of the integrals has been considered.
More precisely, the densities have been approximitated by means of 2-dimension \emph{Chebyshev polynomials}\,, which have several good properties:
\begin{itemize}
\item The approximation grows quickly with the degree of the polynomial, without oscilliation phenomena,
\item The Chebyshev transform is quite related to the fourier transform, which makes the parameters of the polynoms really quickly computable by means of a Fast Fourier Transform,
\item Integration is easy to compute.
\end{itemize}
In our tests, we have chosen a Chebyshev approximation of degree $128\times 128$\,, which is more than sufficient for an almost exact computation.
\paragraph{Example.}
Two bba $m_1$ and $m_2$ have been constructed by normalizing the following functions $mm_1$ and $mm_2$ defined over $[-1,1]^2$\,:
$$
mm_1\bigl(\grave{x}\cap\acute{y}\bigr)=\exp\bigl(-(x+1)^2-y^2\bigr)
$$
and
$$
mm_2\bigl(\grave{x}\cap\acute{y}\bigr)=\exp\bigl(-x^2-(y-1)^2\bigr)\;.
$$
The fused bba $m_1\oplus m_2$ and the respective belief function $b_1, b_2, b_1\oplus b_2$  have been computed.
This computation has been instantaneous.
All functions have been represented in the figures~\ref{f2k5:fig:mm1} to~\ref{f2k5:fig:b1+b2}.
\paragraph{Interpretation.}
The bba $m_1$ is a density centered around the interval $[-1,0]$\,, while $m_2$ is a density centered around $[0,1]$\,.
This explains why the belief $b_1$ increases faster from the interval $[-1,-1]$ to $[-1,1]$ than from the interval $[1,1]$ to $[-1,1]$\,.
And this property is of course inverted for~$b_2$\,.\\[3pt]
A comparison of the fused bba $m_1\oplus m_2$ with the initial bbas $m_1$ and $m_2$ makes apparent a global forward move of the density.
This just means that the fused bba is put on intervals with less imprecision, and possibly on some intervals with negative width (\emph{ie.} associated with a degree of contradiction).
Of course there is nothing surprising here, since information fusion will reduce imprecision and produce some contradiction!
It is also noticed that the fused bba is centered around the interval $[0,0]$\,.
This result matches perfectly the fact that $m_1$ and $m_2$\,, and their related sensors, put more belief respectively over the interval $[-1,0]$ and the interval $[0,1]$\,; and of course $[-1,0]\cap[0,1]=[0,0]$\,.
\begin{figure}
\centering
\scalebox{.75}{\input{Test_plotData.mm2.tex}}
\caption{Non normalized bba $mm_2$}
\label{f2k5:fig:mm2}
\end{figure}
\begin{figure}
\centering
\scalebox{.75}{\input{Test_plotData.m1.tex}}
\caption{Basic belief assignment $m_1$}
\label{f2k5:fig:m1}
\end{figure}
\begin{figure}
\centering
\scalebox{.75}{\input{Test_plotData.m2.tex}}
\caption{Basic belief assignment $m_2$}
\label{f2k5:fig:m2}
\end{figure}
\begin{figure}
\centering
\scalebox{.75}{\input{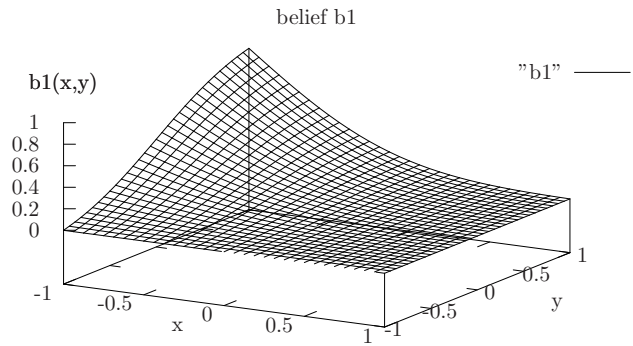}}
\caption{Belief function $b_1$}
\label{f2k5:fig:b1}
\end{figure}
\begin{figure}
\centering
\scalebox{.75}{\input{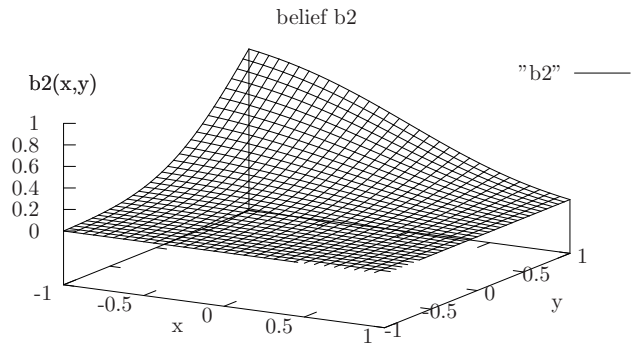}}
\caption{Belief function $b_2$}
\label{f2k5:fig:b2}
\end{figure}
\begin{figure}
\centering
\scalebox{.75}{\input{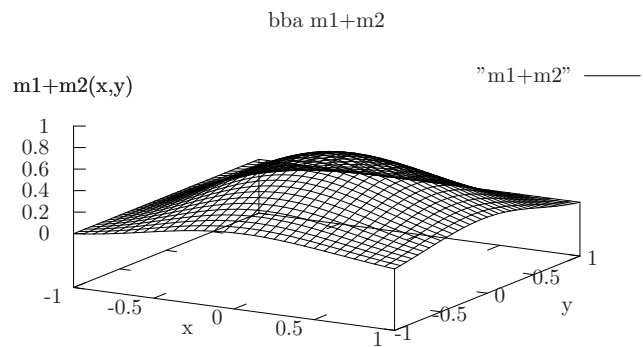}}
\caption{Fused bba $m_1\oplus m_2$}
\label{f2k5:fig:m1+m2}
\end{figure}
\begin{figure}
\centering
\scalebox{.75}{\input{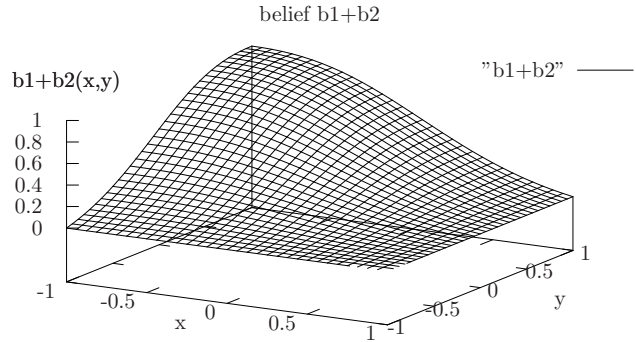}}
\caption{Fused bba $b_1\oplus b_2$}
\label{f2k5:fig:b1+b2}
\end{figure}
%
%
%
\section{Conclusion}
\label{f2k5:sectionConclusion}
A problem of continuous information fusion has been investigated and solved in the DSmT paradigm.
The conceived method is based on the generalization of the notion of hyperpower set.
It is versatile and is able to specify the typical various degrees of contradiction of a DSm model.
It has been implemented efficiently for a bounded continuous information.
The work is still prospective, but applications should be done in the future on localization problems.
At this time, the concept is restricted to one-dimension informations.
However, works are now accomplished in order to extend the method to multiple-dimensions domains.
%
%
\appendix
\section{Biography}
{\bf FR\'ED\'ERIC DAMBREVILLE}
Frederic Dambreville studied mathematics, logic, signal and image processing. He received the Ph.D.
degree in signal processing and optimization, from the university of Rennes, France, in 2001.
He enjoyed a stay in California (U.S.A.) and worked as a postdoctorate in the Naval Postgraduate School at Monterey in 2001/2002.
In 2002, he joined the department image, perception and robotic of the CTA laboratory (Delegation Generale pour l'Armement), France.
His main interests are in optimization, optimal planning, game theory, simulation methods, data\&sensor fusion, Markov models\&Bayesian networks, Logic\&Conditional logic.
His most recent works are about rare event simulation (\emph{e.g.} Cross-Entropy optimization), optimal decision with partial observation, hierachical system optimization, Sche\-duling, modal\&Bayesian logic, and DSmT.
\\[5pt]
D\'el\'egation G\'en\'erale pour l'Armement, DGA/CEP/GIP/SRO\\
\noindent
16 Bis, Avenue Prieur de la C\^ote d'Or\\
Arcueil, F 94114, France\\
Web: {\tt http://www.FredericDambreville.com}\\
Email: {\tt http://email.FredericDambreville.com}\\[5pt]
%

%
\end{document}